\documentclass{amsart}
\usepackage{amsmath}
\usepackage{graphicx}

\newtheorem{theorem}{Theorem}

\newtheorem{corollary}[theorem]{Corollary}

\newtheorem{lemma}[theorem]{Lemma}

\newcommand{\sg}{\sigma}

\begin{document}
\title{Cohomology of topological graphs and Cuntz-Pimsner algebras}
\date{\today}
\author{V.~Deaconu}
\address{Dept of Math/084\\
Univ of Nevada\\
Reno, NV 89557}
\email{vdeaconu@unr.edu}
\thanks{Research of the first two authors was supported by \textsc{nsf} grant
DMS-9706982.}
\author{A.~Kumjian}
\address{Dept of Math/084\\
Univ of Nevada\\
Reno, NV 89557}
\email{alex@unr.edu}
\author{P. Muhly}
\address{Dept of Math\\
Univ of Iowa\\
Iowa City, IA 52242}
\email{muhly@math.uiowa.edu}
\thanks{Research of the third author was supported by \textsc{nsf} grant DMS
97-06713.}

\begin{abstract}
We compute the sheaf cohomology of a groupoid built from a local
homeomorphism of a locally compact space $X$. In particular, we identify the
twists over this groupoid, and its Brauer group. Our calculations refine
those made by Kumjian, Muhly, Renault and Williams in the case $X$ is the
path space of a graph, and the local homeomorphism is the shift. We also
show how the C*-algebra of a twist may be identified with the Cuntz-Pimsner
algebra constructed from a certain C*-correspondence.
\end{abstract}

\maketitle


\section{Introduction\protect\bigskip}

Let $X$ be a second countable, locally compact, Hausdorff space and let $%
\sigma :X\rightarrow X$ be a local homeomorphism (not necessarily
surjective). The pair $(X,\sigma )$ will be fixed throughout this note. In 
\cite[Theorem 1]{d1}, the first author showed how to build an $r$-discrete
groupoid with Haar system from $(X,\sigma )$ as follows: 
\begin{equation*}
\Gamma =\Gamma (X,\sigma )=\{(x,m,y)\in X\times \mathbf{Z}\times X~|%
\mbox{ there are $k, l \geq 0$ such that $\sg^k (x) = \sg^l (y)$ and
$k - l = m$ }\}.
\end{equation*}
The groupoid operations are given by the formulae: 
\begin{gather*}
r(x,m,y)=x,\quad s(x,m,y)=y,\quad (x,m,y)(y,n,z)=(x,m+n,z),\quad \\
(x,m,y)^{-1}=(y,-m,x).
\end{gather*}
The Haar system is given by the counting measures on the sets $r^{-1}(x)$, $%
x\in X$. (See \cite{A-D} and \cite{AR} also.) As noted in \cite{d1},
groupoids of the form $\Gamma $ generalize transformation groupoids
associated with homeomorphisms. Indeed, if $\sigma $ is a homeomorphism,
then $\Gamma $ is its transformation group groupoid. For another example to
keep in mind in this note, suppose that $E$ is a graph with no sinks as in 
\cite{kpr}, suppose that $X=E^{\infty }$, the infinite path space, and
suppose that $\sigma :X\rightarrow X$ is the unilateral shift, 
\begin{equation*}
\sigma (x_{1},x_{2},x_{3},\dots )=(x_{2},x_{3},\dots )\text{.}
\end{equation*}
Then $\sigma $ is a local homeomorphism, and $\Gamma =\mathcal{G}_{E}$ - the
groupoid studied in \cite{kpr}.

The aim of this note is twofold. First we use the long exact sequence of 
\cite[3.7]{k2} to compute the sheaf cohomology of $\Gamma$. This computation
allows us to identify explicitly all (circle) twists over $\Gamma$ in the
sense of \cite{k1}. These are extensions $\Lambda$ of $\Gamma$ by the
groupoid $X\times\mathbf{T}$. Roughly speaking, twists over groupoids
generalize $2$-cocycles over groups and the (restricted) groupoid $C^{\ast}$%
-algebra of a twist (to be defined below) is the groupoid analogue of the $%
C^{\ast}$-algebra of a group twisted by a $2$-cocycle. Our principal result
in this direction, Theorem \ref{cohomiso}, which is proved in the next
section, asserts that for any sheaf of abelian groups $A$ on which $\Gamma$
acts, $Z_{\Gamma}^{n}(A)$ is naturally isomorphic to $H^{n}(X,A)$, where $%
Z_{\Gamma}^{n}$ is the $n^{\mathrm{th}}$ right derived functor of the
cocycle functor and $H^{n}(X,A)$ is the usual sheaf cohomology of $X$ with
values in $A$. Further, our computation of the sheaf cohomology of $\Gamma$
allows us to refine the calculations made in \cite[Proposition 11.8]{kmrw}
that show that the Brauer group of $\Gamma$ vanishes when $X$ is the path
space of a graph (with no sinks). In Section 4, we give some other examples
that illustrate this in settings that are of current interest in operator
algebra.

Our second objective is to use the calculations of Section 2 to show that
for each twist $\Lambda $ over $\Gamma $, the (restricted) groupoid $C^{\ast
}$-algebra $C^{\ast }(\Gamma ;\Lambda )$ is naturally isomorphic to a
Cuntz-Pimsner algebra (see \cite{p}) constructed from the data used to build 
$\Lambda $. A bit more explicitly, first note that $\Lambda $ is a bona fide
groupoid with Haar system in its own right and we may therefore form its $%
C^{\ast }$-algebra, $C^{\ast }(\Lambda )$. It is a completion of the space
of continuous, complex-valued, compactly supported functions on $\Lambda $, $%
C_{c}(\Lambda )$. The circle $\mathbf{T}$ acts on $\Lambda $ in the obvious
way, and the restricted groupoid $C^{\ast }$-algebra of $\Lambda $, $C^{\ast
}(\Gamma ;\Lambda )$, is defined to be the closure in $C^{\ast }(\Lambda )$
of $\{f\in C_{c}(\Lambda )\mid f(z\lambda )=zf(\lambda ),\;z\in \mathbf{T}\}$%
. For the notion of Cuntz-Pimsner algebras, recall from \cite{p} and \cite
{MS1} that a $C^{\ast }$\emph{-correspondence} $\mathcal{E}$ over a $C^{\ast
}$-algebra $A$ is a (right) Hilbert $C^{\ast }$-module $\mathcal{E}$ over $A$
that is endowed with a $C^{\ast }$-representation $\varphi $ of $A$ into the
space of continuous, adjointable module maps on $\mathcal{E}$, $\mathcal{L}(%
\mathcal{E})$. We shall assume that our correspondences are faithful,
meaning that $\varphi $ is injective. However, we do not assume that they
are full, meaning that the closed span of $\langle \mathcal{E},\mathcal{E}%
\rangle $, which is an ideal in $A$, is, in fact, all of $A$. Also, in order
to lighten the notation, we shall usually not mention $\varphi $ unless it
helps to clarify some point. Given a correspondence $\mathcal{E}$ one can
build a $C^{\ast }$-algebra from $\mathcal{E}$, denoted $\mathcal{O}_{%
\mathcal{E}}$, that is a simultaneous generalization of a crossed product of 
$A$ determined by an automorphism and of a Cuntz-Krieger algebra. In \cite{p}%
, $\mathcal{O}_{\mathcal{E}}$ is called the $C^{\ast }$-algebra of $\mathcal{%
E}$, while in \cite{MS1} and elsewhere, $\mathcal{O}_{\mathcal{E}}$ is
called the \emph{Cuntz-Pimsner algebra} associated with $\mathcal{E}$. The
faithfulness of $\varphi $ guarantees that $\mathcal{O}_{\mathcal{E}}$ is
non-zero. In \cite{p}, Pimsner usually assumes that his modules $\mathcal{E}$
are full. However, this is not necessary for our purposes and we choose to
use what he calls the \emph{augmented} algebra in \cite[Remark 1.2.3]{p}
instead of his $\mathcal{O}_{\mathcal{E}}$. We shall not need anything about
the actual construction of $\mathcal{O}_{\mathcal{E}}$ in this note and we
shall only use a few of the properties of these algebras. We shall therefore
refer the reader to \cite{p} and \cite{MS1} for most details.

If $X$ is compact it was shown by the first author in \cite[3.3]{d2} that $%
C^{\ast}(\Gamma)$ may be identified with the Cuntz-Pimsner algebra $\mathcal{%
O}_{\mathcal{H}}$ where $\mathcal{H}$ is the $C^{\ast}$-correspondence over $%
C(X)$, denoted $\ell^{2}(\sigma)$, naturally associated to $\sigma$. His
analysis works even when $X$ is locally compact. In this case, $%
\ell^{2}(\sigma)$ is defined to be the completion of the pre-Hilbert $%
C^{\ast}$-module $C_{c}(X)$ over $C_{0}(X)$ defined by the formulae: 
\begin{equation*}
\xi\cdot f(x)=\xi(x)f(\sigma(x))
\end{equation*}
and 
\begin{equation*}
\langle\xi,\eta\rangle(x)=\sum_{\sigma(y)=x}\overline{\xi(y)}\eta(y)\text{,}
\end{equation*}
$\xi, \eta \in C_{c}(X)$, $f\in C_{0}(X)$. The fact that $\sigma$ is a local
homeomorphism, coupled with the fact that $\xi$ and $\eta$ both lie in $%
C_{c}(X)$ guarantee that the sum defining the inner product is finite. The
algebra $C_{0}(X)$ acts on $\ell^{2}(\sigma)$ to the left via the formula 
\begin{equation*}
f\cdot\xi(x)=f(x)\xi(x)\text{,}
\end{equation*}
and with respect to this left action $\ell^{2}(\sigma)$ becomes a $C^{\ast}$%
-correspondence over $C_{0}(X)$.

As we shall show in Section 3, twists over $\Gamma $ are naturally
associated to line bundles over $X$. If $T$ is such a line bundle, then we
may build the associated twist $\Lambda _{T}$ over $\Gamma $ and we may
``twist'' the correspondence $\ell ^{2}(\sigma )$ by $T$ to obtain a new $%
C^{\ast }$-correspondence $\mathcal{H}_{T}$ over $C_{0}(X)$. We shall show
in Theorem \ref{CPident} that there is a natural isomorphism between $%
C^{\ast }(\Gamma ;\Lambda _{T})$ and the Cuntz-Pimsner algebra $\mathcal{O}_{%
\mathcal{H}_{T}}$. If the map $\sigma $ is a homeomorphism, then the
correspondences $\ell ^{2}(\sigma )$ and $\mathcal{H}_{T}$ are
``invertible'', meaning that they are Hilbert bimodules in the sense used by
Abadie, Eilers, and Exel in \cite{aee} and by others. In this event, the
presentation in \cite{aee} shows how to realize $\mathcal{O}_{\mathcal{H}%
_{T}}$ as a Fell bundle over $\mathbf{Z}$. Our analysis in the general
setting provides a way of thinking about $\mathcal{O}_{\mathcal{H}_{T}}$ in
terms of a twist over $\Gamma $ - a kind of a Fell bundle over $\Gamma $.

One benefit of our isomorphism theorem, Theorem \ref{CPident}, that is under
development, is that we will be able to apply the results from \cite{MS2}
and \cite{fmr} to give conditions implying that $C^{\ast }(\Gamma ;\Lambda )$
is simple. We also will be able to apply technology developed in \cite{d1}
to calculate the $K$-theory of $C^{\ast }(\Gamma ;\Lambda )$.

\section{Cohomology Calculations}

A $\Gamma$-sheaf is simply a sheaf over the unit space $X$ of $\Gamma$ on
which $\Gamma$ acts. We shall view a sheaf $A$ over $X$ both as an \'{e}tale
space over $X$ with abelian group fibers and as a functor from the category
of open subsets of $X$ to the category of abelian groups satisfying the
usual relations. For an \emph{arbitrary }$r$-discrete groupoid $\Gamma$ and $%
\Gamma$-sheaf\ $A$ one has the following long exact sequence (by \cite[3.7]
{k2}): 
\begin{gather}
0\rightarrow H^{0}(\Gamma,A)\rightarrow H^{0}(\Gamma^{0},A)\overset{d}{%
\longrightarrow}Z_{\Gamma}^{0}(A)\rightarrow H^{1}(\Gamma,A)\rightarrow
\cdots  \notag \\
\rightarrow H^{n-1}(\Gamma^{0},A)\overset{d}{\longrightarrow}%
Z_{\Gamma}^{n-1}(A)\rightarrow H^{n}(\Gamma,A)\rightarrow
H^{n}(\Gamma^{0},A)\rightarrow Z_{\Gamma}^{n}(A)\rightarrow\cdots
\label{1st les}
\end{gather}
where $H^{n}(\Gamma,A)$ denotes the $n^{\mathrm{th}}$ equivariant cohomology
of $\Gamma$ with coefficients in\ $A$ (cf.~\cite{g}), $H^{n}(\Gamma^{0},A)$
denotes the usual sheaf cohomology of the unit space $\Gamma^{0}$ with
coefficients in $A$ (we use the same symbol for a $\Gamma$-sheaf\ and its
underlying sheaf) and $Z_{\Gamma}^{n}$ denotes the $n^{\mathrm{th}}$ right
derived functor of the cocycle functor. The cocycle functor $Z_{\Gamma}:%
\mbox{\sf Ab}%
%
(\Gamma)\rightarrow%
\mbox{\sf Ab}%
%
$, where $%
\mbox{\sf Ab}%
%
(\Gamma)$ is the category of $\Gamma$-sheaves and $%
\mbox{\sf Ab}%
%
$ is the category of abelian groups, is defined as follows: Given a $\Gamma $%
-sheaf\ $A$, the abelian group $Z_{\Gamma}(A)$ consists of all continuous
functions $f:\Gamma\rightarrow A$ such that $f(\gamma)\in A_{r(\gamma)}$
(i.e.\ $f$ is a continuous section of $r^{\ast}(A)$) and $%
f(\gamma_{1}\gamma_{2})=f(\gamma_{1})+\gamma_{1}f(\gamma_{2})$ for all $%
(\gamma_{1},\gamma_{2})\in\Gamma^{2}$. Thus, $Z_{\Gamma}(A)$ is the usual
group of one-cocycles or crossed homomorphisms with values in the bundle $A$.

When $\Gamma=\Gamma(X,\sigma)$, we have $\Gamma^{0}=X$ (under the
identification $(x,0,x)\mapsto x$) and we shall show that there is an
isomorphism $Z_{\Gamma}^{n}(A)\simeq H^{n}(X,A)$ for any $\Gamma$-sheaf\ $A$.

Note that $H^{0}(X,A)=S(A)$, the group of continuous sections of $A$, and
that $H^{n}(X,\cdot)$ is the $n^{\mathrm{th}}$ right derived functor of $S$
(when it is regarded as a functor from the category of sheaves of abelian
groups over $X$, $%
\mbox{\sf Ab}%
%
(X)$ to $%
\mbox{\sf Ab}%
%
$.) Our first goal is to show that the functors $S$ and $Z_{\Gamma}$ are
naturally isomorphic. Define a map $\varphi_{A}:Z_{\Gamma}(A)\rightarrow
S(A) $ by $\varphi_{A}(f)(x)=f(x,1,\sigma(x))$.

\begin{lemma}
\label{naturaliso}The map $\varphi _{A}$ defines a natural isomorphism
between $Z_{\Gamma }$ and $S$.
\end{lemma}

\begin{proof}
It is easy to check that $\varphi _{A}$ defines a natural transformation
between the functors $Z_{\Gamma }$ and $S$ and that $\varphi _{A}$ is a
homomorphism. It remains to show that it is bijective. The map $x\mapsto
(x,1,\sigma (x))$ defines a homeomorphism from $X$ onto a clopen subset of $%
\Gamma ,$ which we denote by $X^{\prime }$. The bijectivity of $\varphi _{A}$
is equivalent to the assertion that every continuous section on $X^{\prime }$
has a unique continuous extension to $\Gamma $ that satisfies the cocycle
property. But this is straightforward: For $\gamma =(x,k-l,y)\in \Gamma $
(with $\sigma ^{k}(x)=\sigma ^{l}(y)$) one has the factorization $\gamma
=\xi _{1}\cdots \xi _{k}\eta _{l}^{-1}\cdots \eta _{1}^{-1}$ where $\xi
_{i}=(\sigma ^{i-1}(x),1,\sigma ^{i}(x))$ for $i=1,\dots ,k$ and $\eta
_{j}=(\sigma ^{j-1}(y),1,\sigma ^{j}(y))$ for $j=1,\dots ,l$. Note that $\xi
_{i},\eta _{j}\in X^{\prime }$ and that the extension is uniquely determined
by the cocycle property $f(\gamma _{1}\gamma _{2}\cdots \gamma
_{n})=f(\gamma _{1})+\gamma _{1}f(\gamma _{2})+\cdots +\gamma _{1}\gamma
_{2}\cdots \gamma _{n-1}f(\gamma _{n})$ and the fact that $f(\gamma
^{-1})=-\gamma ^{-1}f(\gamma )$. It follows that $\varphi _{A}$ is indeed
bijective and so defines a natural isomorphism between the functors $%
Z_{\Gamma }$ and $S$.
\end{proof}

\begin{theorem}
\textbf{\label{cohomiso}} The map $\varphi _{A}$ induces an isomorphism
between $Z_{\Gamma }^{n}(A)$ and $H^{n}(X,A)$.
\end{theorem}

\begin{proof}
Given a $\Gamma $-sheaf\ $B$, there is a $\Gamma $-sheaf\ $Q(B)$ (see 
\cite[1.6]{k2}) which is flabby as a sheaf over $X$ (a sheaf is flabby or
flasque if any continuous section defined on an open subset of $X$ may be
extended continuously to all of $X$). Hence, for any $\Gamma $-sheaf\ $A$
there is an injective resolution (in $%
\mbox{\sf Ab}%
%
(\Gamma )$) $A\rightarrow Q^{0}\rightarrow Q^{1}\rightarrow Q^{2}\rightarrow
\dots $, which is flabby when regarded as a resolution in $%
\mbox{\sf Ab}%
%
(X)$. By applying the functors $Z_{\Gamma }$ and $S$ to the complex $Q^{*}$
and invoking Lemma \ref{naturaliso}, one obtains the diagram: 
\begin{equation*}
\begin{array}{ccccccc}
Z_{\Gamma }(Q^{0}) & \rightarrow & Z_{\Gamma }(Q^{1}) & \rightarrow & 
Z_{\Gamma }(Q^{2}) & \rightarrow & \cdots \\ 
\downarrow &  & \downarrow &  & \downarrow &  &  \\ 
S(Q^{0}) & \rightarrow & S(Q^{1}) & \rightarrow & S(Q^{2}) & \rightarrow & 
\cdots
\end{array}
\end{equation*}
in which the vertical arrows are isomorphisms $\varphi _{Q^{i}}:Z_{\Gamma
}(Q^{i})\rightarrow S(Q^{i})$ and the diagram commutes because $\varphi _{A}$
is a natural transformation. Then $Z_{\Gamma }^{n}(A)$ is the cohomology of
the complex corresponding to the first row. On the other hand, $H^{n}(X,A)$
is the cohomology of the second row since $H^{n}(X,F)=0$ if $F$ is flabby
and $n>0$ (see \cite[3.15]{t}, \cite[II.3.5]{i}).
\end{proof}

\begin{corollary}
\label{les}The long exact sequence (\ref{1st les}) induces the long exact
sequence 
\begin{gather}
0\rightarrow H^{0}(\Gamma ,A)\rightarrow H^{0}(X,A)\overset{1-\sigma ^{\ast }%
}{\longrightarrow }H^{0}(X,A)\rightarrow H^{1}(\Gamma ,A)\rightarrow \cdots 
\notag \\
\rightarrow H^{n-1}(X,A)\overset{1-\sigma ^{\ast }}{\longrightarrow }%
H^{n-1}(X,A)\rightarrow H^{n}(\Gamma ,A)\rightarrow  \notag \\
H^{n}(X,A)\overset{1-\sigma ^{\ast }}{\longrightarrow }H^{n}(X,A)\rightarrow
\cdots  \label{2nd les}
\end{gather}
where $\sigma ^{\ast }$ is the map on cohomology induced by the local
homeomorphism $\sigma :X\rightarrow X$.
\end{corollary}

\begin{proof}
We write $\sigma ^{\ast }$ also for the pull back of sheaves. That is, $%
\sigma ^{\ast }(A)$ is the pull back sheaf on $X$ induced by $A$ and $\sigma 
$. It is isomorphic to $A$ since $A$ is a $\Gamma $-sheaf. Conversely, given
a sheaf $B$ over $X$ together with an isomorphism $B\simeq \sigma ^{\ast
}(B) $ one may endow $B$ with the structure of a $\Gamma $-sheaf in a
natural way. For the $\Gamma $-sheaf\ $A$, the long exact sequence (\ref{1st
les}) arises from a map $d:S(A)\rightarrow Z_{\Gamma }$ defined by the
equation $d(f)(\gamma )=f(r(\gamma ))-\gamma f(s(\gamma ))$. Since the
composition $\varphi _{A}d:S(A)\rightarrow S(A)$ is given by the equation 
\begin{align*}
\varphi _{A}d(f)(x)& =f(r(x,1,\sigma (x)))-(x,1,\sigma (x))f(s(x,1,\sigma
(x))) \\
& =f(x)-(x,1,\sigma (x))f(\sigma (x))\text{,}
\end{align*}
it follows that $\varphi _{A}d=1-\sigma ^{\ast }$. With this observation, we
see that the long exact sequence (\ref{1st les}) may be rewritten as the
long exact sequence (\ref{2nd les}).
\end{proof}

In the case of most interest to us, $A$ is the sheaf $\mathcal{S}$ of germs
of continuous circle-valued functions on $X$. Since $\Gamma $ is $r$%
-discrete, elements in $\Gamma $ may be viewed as germs of local
homeomorphisms of $X$. Hence there is a canonical action of $\Gamma $ on $%
\mathcal{S}$ given by composition of germs. There is an extension of sheaves 
\begin{equation*}
0\rightarrow \mathbf{Z}\rightarrow \mathcal{R}\rightarrow \mathcal{S}%
\rightarrow 0
\end{equation*}
where $\mathcal{R}$ is the sheaf of germs of continuous real-valued
functions on $X$ and $\mathbf{Z}$ is the constant sheaf of integers (both
endowed with canonical $\Gamma $-actions). One has the long exact sequence: 
\begin{gather}
0\rightarrow H^{0}(X,\mathbf{Z})\rightarrow H^{0}(X,\mathcal{R})\rightarrow
H^{0}(X,\mathcal{S})\rightarrow H^{1}(X,\mathbf{Z})\rightarrow  \notag \\
H^{1}(X,\mathcal{R})\rightarrow H^{1}(X,\mathcal{S})\rightarrow H^{2}(\Gamma
,\mathbf{Z})\rightarrow \cdots \text{.}
\end{gather}
Since $\mathcal{R}$ is soft, we have $H^{n}(X,\mathcal{R})=0$ for $n>0$ and,
hence, $H^{n}(X,\mathcal{S})\simeq H^{n+1}(X,\mathbf{Z})$, for $n>0$. Since
one also has a short exact sequence of $\Gamma $-sheaves, there is also the
long exact sequence (see \cite[Def. 0.11]{k3}): 
\begin{gather}
0\rightarrow H^{0}(\Gamma ,\mathbf{Z})\rightarrow H^{0}(\Gamma ,\mathcal{R}%
)\rightarrow H^{0}(\Gamma ,\mathcal{S})\rightarrow H^{1}(\Gamma ,\mathbf{Z}%
)\rightarrow  \notag \\
H^{1}(\Gamma ,\mathcal{R})\rightarrow H^{1}(\Gamma ,\mathcal{S})\rightarrow
H^{2}(\Gamma ,\mathbf{Z})\rightarrow \cdots .
\end{gather}
By Corollary \ref{les} and the fact that $\mathcal{R}$ is soft, $%
H^{n}(\Gamma ,\mathcal{R})=0$ for $n>1$; it follows that $H^{n}(\Gamma ,%
\mathcal{S})\simeq H^{n+1}(\Gamma ,\mathbf{Z})$ for all $n>1$. In \cite[11.3]
{kmrw} (see also \cite[4.19]{k2}) the second cohomology group $H^{2}(\Gamma ,%
\mathcal{S})$ was identified with the so-called Brauer group $%
\mbox{Br}%
%
(\Gamma )$ of $\Gamma $. This is the collection of strong Morita equivalence
classes of $\Gamma$-bundles of elementary $C^{\ast }$-algebras satisfying
Fell's condition. These equivalence classes form a group under tensor
product that generalizes the Brauer group of finite dimensional, central
simple algebras over a field. These facts together with Corollary \ref{les}
yields the following:

\begin{corollary}
\label{les bis} We have $\mathrm{Br}(\Gamma )\simeq H^{2}(\Gamma ,\mathcal{S}%
)\simeq H^{3}(\Gamma ,\mathbf{Z})$. Hence, in the notation of Corollary \ref
{les}, one has the following exact sequence: 
\begin{equation}
\begin{array}{ccccccccc}
H^{2}(X,\mathbf{Z}) & \overset{1-\sigma ^{\ast }}{\longrightarrow } & 
H^{2}(X,\mathbf{Z}) & \rightarrow & \mathrm{Br}(\Gamma )\rightarrow & 
H^{3}(X,\mathbf{Z}) & \overset{1-\sigma ^{\ast }}{\longrightarrow } & 
H^{3}(X,\mathbf{Z}) & 
\end{array}
\label{les brauer}
\end{equation}
\end{corollary}

Note that this includes the fact that $%
\mbox{Br}%
%
(\Gamma )=0$ when $X$ is the path space of a graph and $\sigma$ is the
unilateral shift \cite[Proposition 11.8]{kmrw}. \bigskip

\section{Cuntz-Pimsner Algebras}

Suppose $\Gamma$ is a general $r$-discrete groupoid and that $A$ is a $%
\Gamma $-sheaf. Then a \emph{twist} by $A$ over $\Gamma$ is an $r$-discrete
groupoid $\Sigma$, with $\Sigma^{0}=\Gamma^{0}$, together with two groupoid
homomorphisms, $j$ and $\pi$, so that 
\begin{equation*}
A\underset{j}{\longrightarrow}\Sigma\underset{\pi}{\longrightarrow}\Gamma%
\text{,}
\end{equation*}
with $j$ injective and $\pi$ surjective, $\pi^{-1}(\Gamma^{0})=j(A)$, and so
that $\sigma j(a)\sigma^{-1}=j(\pi(\sigma)a)$ for all $\sigma\in\Sigma$ and
all $a\in A_{s(\sigma)}$. Thus a twist by $A$ over $\Gamma$ is simply an
extension of $\Gamma$ by $A$. Two twists are called isomorphic in case they
are isomorphic as extensions in the usual sense. The isomorphism classes of
twists by $A$ over $\Gamma$ becomes an abelian group under Baer sum that is
denoted $T_{\Gamma}(A)$.

In the special case when $A=\mathcal{S}$, $T_{\Gamma }(\mathcal{S})$ is also
written $%
\mbox{Tw}%
%
(\Gamma )$. Further, as we mentioned in the introduction, a twist $\Sigma $
by $\mathcal{S}$ over $\Gamma $ may be viewed as a principal circle bundle
over $\Gamma $ where the circle action is compatible with the groupoid
actions. For the details on the theory of twists, see \cite{k1} and 
\cite[Section 2]{k2}. Note, however, that in \cite{k1}, twists are
restricted to extensions by $\mathcal{S}$ of \emph{principal }groupoids. The
restriction to principal groupoids is not necessary for our purposes and
indeed, in general, $\Gamma (X,\sigma )$ is not principal. In Corollary 3.4
of \cite{k2}, it is proved that $T_{\Gamma }(A)$ is naturally isomorphic to $%
Z_{\Gamma }^{1}(A)$ for any $\Gamma $-sheaf $A$. Hence, when $\Gamma =\Gamma
(X,\sigma )$ and $A=\mathcal{S}$ we conclude from Theorem \ref{cohomiso}
that 
\begin{equation*}
H^{1}(X,\mathcal{S})\simeq Z_{\Gamma }^{1}(\mathcal{S})\simeq T_{\Gamma }(%
\mathcal{S})\simeq 
\mbox{Tw}%
%
(\Gamma ).
\end{equation*}
Our objective in this section is to construct this isomorphism directly and
then to show that the twisted groupoid $C^{\ast }$-algebra is a
Cuntz-Pimsner algebra.

Recall that $H^{1}(X,\mathcal{S})$ may be identified with the group of
isomorphism classes of principal circle bundles over $X$. We want to see how
to pass between circle bundles over $X$ to twists -- i.e., certain circle
bundles over $\Gamma$. To this end we shall write $j$ for the map that sends 
$x\in X$ to $(x,1,\sigma(x))$ in $\Gamma$ (see Lemma \ref{naturaliso}); then 
$j$ induces a homeomorphism from $X$ to its image, the clopen subset $%
X^{\prime}:=\{(x,1,\sigma(x))~|~x\in X\}$. (This $j$ should not be confused
with the $j$ discussed above in the general theory of twists, which will
never be mentioned again.) Given a twist $\Lambda$ over $\Gamma$, viewed as
a principal circle bundle over $\Gamma$, we may pull the circle bundle back
to $X$, via $j$, to obtain a principal circle bundle over $X$. In symbols, $%
\Lambda\mapsto j^{\ast}(\Lambda)$. 
We are thus led to

\begin{theorem}
\label{twistiso}The map $\Lambda \rightarrow j^{\ast }(\Lambda )$ implements
an isomorphism from $%
\mbox{Tw}%
%
(\Gamma )$ onto $H^{1}(X,\mathcal{S})$ viewed as isomorphism classes of
principal circle bundles over $X$.
\end{theorem}

\begin{proof}
Our comments prior to the statement of the theorem together with a moment's
reflection reveal that the map $\Lambda \rightarrow j^{\ast }(\Lambda )$
induces a homomorphism from $\mbox{Tw}(\Gamma )$ into $H^{1}(X,\mathcal{S})$%
. It remains to show that the map is bijective. We first show the
surjectivity, i.e., how to construct a twist $\Lambda _{T}$ from a principal
circle bundle $T$ so that $T\simeq j^{\ast }(\Lambda _{T})$. So let the
principal circle bundle $T$ over $X$ be given, write $p:T\rightarrow X$ for
the quotient map, and for $k,l\geq 0$ set 
\begin{equation*}
X_{k,l}=\{(x,k-l,y)~|~\sigma ^{k}(x)=\sigma ^{l}(y)\}\subset \Gamma .
\end{equation*}
It is straightforward to verify the following assertions:

\begin{enumerate}
\item  $X_{k,l}$ is a clopen subset of $\Gamma $,

\item  $X_{k,l}\subset X_{k+1,l+1},$ and

\item  $c^{-1}(k-l)=\cup _{j=0}^{\infty }X_{k+j,l+j}$, where $c$ is the
position cocycle: $c(x,k-l,y)=k-l$.
\end{enumerate}

Our strategy is to ``extend'' $T$ to $X_{k,l}$ in such a way that the
extensions to $X_{k,l}$ and to $X_{k+1,l+1}$ are compatible. This will give
bundles over the disjoint sets, $c^{-1}(n)$, $n\in \mathbf{Z}$, which may
then be pieced together in the obvious way. Given circle bundles $T_{1}$ and 
$T_{2}$ over $X_{1}$ and $X_{2}$ one may form the ``product'' circle bundle
over $X_{1}\times X_{2}$, $T_{1}\star T_{2}=T_{1}\times T_{2}/\sim $, where $%
(zt_{1},t_{2})\sim (t_{1},zt_{2})$ for $z\in \mathbf{T}$. We let $\overline{T%
}$ denote the conjugate circle bundle; there is a fiber preserving
homeomorphism $T\rightarrow \overline{T}$, written $t\mapsto \overline{t}$,
such that $\overline{zt}=\overline{z}\overline{t}$ for all $z\in \mathbf{T}$%
. Note that the pull-back of $T\star \overline{T}$ along the diagonal is
canonically isomorphic to the trivial circle bundle $X\times \mathbf{T}$.
Indeed, given $t_{1},t_{2}\in T$ with $p(t_{1})=p(t_{2})$, there is a unique 
$z\in \mathbf{T}$ so that $t_{1}=zt_{2}$; write $z=t_{1}\overline{t_{2}}$.
The desired isomorphism is then given by $(t_{1},t_{2})\mapsto
(p(t_{1}),t_{1}\overline{t_{2}})$ for all $t_{1},t_{2}\in T$ with $%
p(t_{1})=p(t_{2})$ (note that this is well-defined). By a slight abuse of
notation let $T^{k}$ denote the circle bundle $T\star \cdots \star T$ ($k$%
-factors) over $X^{k}$. Observe that there is a natural embedding $\iota
_{k,l}:X_{k,l}\rightarrow X^{k}\times X^{l}$ given by the formula 
\begin{equation*}
\iota _{k,l}(x,m,y)=(x,\sigma (x),\dots ,\sigma ^{k-1}(x),\sigma
^{l-1}(y),\dots ,\sigma (y),y).
\end{equation*}
Set $T_{k,l}=\iota _{k,l}^{\ast }({T}^{k}\star {\overline{T}}^{l})$ (note
that $T_{0,0}$ is the trivial circle bundle over $X_{0,0}=X$). One verifies
that the restriction of $T_{k+1,l+1}$ to $X_{k,l}$ is isomorphic to $T_{k,l}$%
: if $(u_{1},\dots ,u_{k+1},\overline{v_{l+1}},\dots ,\overline{v_{1}})\in
T_{k+1,l+1}$ lies in the fiber over $(x,k-l,y)\in X_{k,l}$, then $%
p(u_{k+1})=\sigma ^{k}(x)=\sigma ^{l}(y)=p(v_{l+1})$ and the desired
isomorphism is given by 
\begin{equation*}
(u_{1},\dots ,u_{k+1},\overline{v_{l+1}},\dots ,\overline{v_{1}})\mapsto
u_{k+1}\overline{v_{l+1}}(u_{1},\dots ,u_{k},\overline{v_{l}},\dots ,%
\overline{v_{1}}).
\end{equation*}
The desired twist $\Lambda _{T}$ is obtained by piecing together the circle
bundles $T_{k,l}$. We also denote the quotient map from $\Lambda _{T}$ to $%
\Gamma $ by $p$. The source and range maps on $\Lambda _{T}$ are defined via 
$p$, i.e.\ $s(\lambda )=s(p(\lambda ))$ and $r(\lambda )=r(p(\lambda ))$ for
all $\lambda \in \Lambda _{T}$. Suppose that $\lambda ,\mu \in \Lambda _{T}$
are composable. Then there are $j,k,l\geq 0$ so that $\lambda \in T_{j,k}$, $%
\mu \in T_{k,l}$ and $x,y,z\in X$ such that $p(\lambda )=(x,j-k,y)$, $p(\mu
)=(y,k-l,z)$. Given $\lambda =(t_{1},\dots ,t_{j},\overline{u_{k}},\dots ,%
\overline{u_{1}})$ and $\mu =(v_{1},\dots ,v_{k},\overline{w_{l}},\dots ,%
\overline{w_{1}})$, define multiplication by the formula, 
\begin{equation*}
\lambda \mu =\left( \prod_{i=1}^{k}v_{i}\overline{u_{i}}\right) (t_{1},\dots
,t_{j},\overline{w_{l}},\dots ,\overline{w_{1}}).
\end{equation*}
(Note that $p(u_{i})=p(v_{i})=\sigma ^{i-1}(y)$, so this formula makes
sense.) It is easy to verify that multiplication is well-defined, continuous
and associative. Finally, the inverse map is defined by 
\begin{equation*}
(u_{1},\dots ,u_{k},\overline{v_{l}},\dots ,\overline{v_{1}}%
)^{-1}=(v_{1},\dots ,v_{l},\overline{u_{k}},\dots ,\overline{u_{1}}).
\end{equation*}
Thus, we see that $\Lambda _{T}$ is a twist over $\Gamma $, and it is
evident that $T=j^{\ast }(\Lambda _{T})$. Hence, the map is surjective.

If the circle bundle is trivial, then any continuous section may be used to
construct a continuous section of the twist (along the above lines) which is
easily seen to be a groupoid homomorphism. Hence, the twist is trivial and
the map is injective.
\end{proof}

Our objective now is to show how to realize $C^{\ast}(\Gamma;\Lambda)$ as a
Cuntz-Pimsner algebra for each twist $\Lambda$ over $\Gamma$. As we have
just seen, each twist $\Lambda$ over $\Gamma$ comes from a unique circle
bundle over $X$. So we begin with these. Given a circle bundle $T$ over $X$,
let $L_{T}$ denote the space of continuous sections of the associated
complex line bundle $T\times_{T}\mathbf{C}$ that vanish at infinity on $T$.
We think of $L_{T}$ as the space of all continuous $\mathbf{C}$-valued
functions $f$ on $T$ such that $f(zt)=zf(t)$ for all $z\in\mathbf{T}$ and
all $t\in T$ and such that $|f|\in C_{0}(X)$. Then, in fact, $L_{T}$ has the
structure of an \emph{imprimitivity} bimodule over $C_{0}(X)$. The action of 
$C_{0}(X)$ is central, i.e., for $f\in C_{0}(X)$ and $\xi\in L_{T}$, $%
f\cdot\xi(t)=\xi\cdot f(t):=\xi(t)f(\dot{t})$; and the inner products are
given by the formulae: ${}_{C_{0}(X)}\langle\xi,\eta\rangle(\dot{t})=\xi(t)%
\overline{\eta}(t)$ and ${}\langle\xi,\eta\rangle_{C_{0}(X)}(\dot{t})=%
\overline{\xi}(t)\eta(t)$. Here $\dot{t}=p(t)$. Note that by the
transformation properties of $\xi$ and $\eta$ these are bona fide $C_{0}(X)$%
-valued inner products. In fact, the most general $C_{0}(X)$-$C_{0}(X)$
equivalence bimodule that fixes the spectrum $X$ is of this form. To say the
same thing differently, the isomorphism classes of these $C_{0}(X)$-$%
C_{0}(X) $ equivalence bimodules form a group under tensor product, a
subgroup of the Picard group of $X$ that is isomorphic to $H^{1}(X,\mathcal{S%
})$ viewed as the isomorphism classes of line bundles over $X$. For these
things, see \cite{PR}.

Next, we ``twist'' $\ell^{2}(\sigma)$ by $L_{T}$. Recall that $%
\ell^{2}(\sigma)$ is a $C^{\ast}$-correspondence over $C_{0}(X)$ and so the
tensor product $\mathcal{H}:=L_{T}\otimes_{C_{0}(X)}\ell^{2}(\sigma)$ makes
sense as a right Hilbert $C^{\ast}$-module over $C_{0}(X)$ (see \cite[5.9]
{ri}). However, since $L_{T}$ is a $C_{0}(X)$-$C_{0}(X)$ imprimitivity
bimodule, the left action of $C_{0}(X)$ on $L_{T}$ passes to one of $%
C_{0}(X) $ on $\mathcal{H}$, making $\mathcal{H}$ a $C^{\ast}$%
-correspondence over $C_{0}(X)$. From the definitions of $L_{T}$ and $%
\ell^{2}(\sigma)$, it is clear that $\mathcal{H}$ may be viewed as the
completion of the compactly supported sections in $L_{T}$ with the following
pre-$C^{\ast}$-correspondence structure: 
\begin{align}
\langle\xi,\eta\rangle_{C_{0}(X)}(x) & =\sum_{\sigma(\dot{t})=x}\overline
{\xi(t)}\eta(t)\text{;}  \label{innerprod} \\
\xi\cdot f(t) & =\xi(t)f(\sigma(\dot{t}))\text{; and}  \notag \\
f\cdot\xi(t) & =f(\dot{t})\xi(t)  \notag
\end{align}
The claim is then that the Pimsner algebra $\mathcal{O}_{\mathcal{H}}$ is
isomorphic to the twisted groupoid $C^{\ast}$-algebra $C^{\ast}(\Gamma
;\Lambda_{T})$ associated to $\Lambda_{T}$, where $\Lambda_{T}$ is the twist
over $\Gamma$ determined by the bundle $T$.

To prove this claim, we need to invoke a result proved in \cite{fmr}. Recall
that if $\mathcal{E}$ is a $C^{\ast}$-correspondence over a $C^{\ast}$%
-algebra $A$, then an\emph{\ (isometric) covariant representation }of $%
\mathcal{E}$ in a $C^{\ast}$-algebra $B$ is a pair $(V,\pi)$ consisting of a 
$\mathbf{C}$-linear map $V:\mathcal{E}\rightarrow B$ and a $C^{\ast}$%
-homomorphism $\pi:A\rightarrow B$ such that the following two conditions
are satisfied:

\begin{enumerate}
\item  $V$ is a bimodule map; i.e., $V(\varphi (a)\xi b)=\pi (a)V(\xi )\pi
(b)$ for all $a,b\in A$ and all $\xi \in \mathcal{E}$.

\item  $V(\xi )^{\ast }V(\eta )=\pi (\langle \xi ,\eta \rangle )$, for all $%
\xi ,\eta \in \mathcal{E}$.
\end{enumerate}

\noindent It can be shown easily that $V$ is bounded and, in fact, $\left\|
V(\xi)\right\| \leq\left\| \xi\right\| $. Furthermore, the map from $%
\mathcal{E}\times\widetilde{\mathcal{E}}$ to $B$ that sends $(\xi,\tilde
{\eta})$ to $V(\xi)V(\eta)^{\ast}$ extends to a $C^{\ast}$-homomorphism $%
\pi^{(1)}$ from $K(\mathcal{E})$, identified with $\mathcal{E}\otimes _{A}%
\widetilde{\mathcal{E}}$, into $B$. (See \cite{p}.) The covariant
representation $(V,\pi)$ is said to satisfy the \emph{Cuntz condition} or to
be a \emph{Cuntz covariant representation} in case $\pi^{(1)}\circ
\varphi(a)=\pi(a)$ for all $a$ in the ideal $J$ in $A$, which is defined to
be $\varphi^{-1}(K(\mathcal{E}))$. It is proved in \cite{p} that an
isometric representation $(V,\pi)$ of $\mathcal{E}$ in a $C^{\ast}$-algebra $%
B$ defines a $C^{\ast}$-representation of the Cuntz-Pimsner algebra $%
\mathcal{O}_{\mathcal{E}}$ in $B$ if and only if $(V,\pi)$ satisfies the
Cuntz condition. The representation of $\mathcal{O}_{\mathcal{E}}$
determined by $(V,\pi)$ is denoted $V\times\pi$ and is called the \emph{%
integrated form} of $(V,\pi)$. Conversely, every representation of the
Cuntz-Pimsner algebra $\mathcal{O}_{\mathcal{E}}$ in a $C^{\ast}$-algebra $B$
is of the form $V\times\pi$ for a (unique) isometric covariant
representation of $\mathcal{E}$ in $B$ that satisfies the Cuntz condition. A
condition for the \emph{faithfulness} of a representation $V\times\pi$ of $%
\mathcal{O}_{\mathcal{E}}$ into $B$, proved in \cite{fmr}, is the following;
it is essential for our analysis.

\begin{lemma}
\label{faithful}Suppose that $(V,\pi )$ is an isometric covariant
representation of $\mathcal{E}$ into a $C^{\ast }$-algebra $B$. Then $%
V\times \pi $ is faithful if and only if $\pi $ is faithful and there is a
(strongly continuous) action $\beta :\mathbf{T}\rightarrow \mbox{\rm Aut}(B)$
such that $\beta _{z}\circ \pi =\pi $ and $\beta _{z}\circ V=zV$ for all $%
z\in \mathbf{T}$.
\end{lemma}

Fix a principal circle bundle $T$ over $X$ and let $\mathcal{H}=\mathcal{H}%
_{T}$ be the $C^{\ast}$-correspondence over $C_{0}(X)$ defined above. Also,
let $j:X\rightarrow \Gamma$ be defined as above by the formula $%
j(x)=(x,1,\sigma(x))$. The bundles $T$ and $\Lambda_{T}$ are related by the
formula 
\begin{equation*}
T=j^{\ast}(\Lambda_{T})\text{.}
\end{equation*}
We define the pair $(V,\pi)$ by the formulae: 
\begin{equation}
\pi(f)(\lambda)=\left\{ 
\begin{array}{cc}
f(p(\lambda)), & \lambda\in\Lambda_{T}|_{\Gamma^{0}} \\ 
0, & \text{otherwise}
\end{array}
\text{,}\right.  \label{v1}
\end{equation}
$f\in C_{0}(X)$, and 
\begin{equation}
V(\xi)(\lambda)=\left\{ 
\begin{array}{cc}
\xi\circ (j^*)^{-1}(\lambda), & \lambda\in\Lambda_{T}|X^{\prime} \\ 
0, & \text{otherwise}
\end{array}
\text{,}\right.  \label{v2}
\end{equation}
$\xi\in\mathcal{H}_{T}$. It is routine to check that $(V,\pi)$ is an
isometric covariant representation of $\mathcal{H}_{T}$ in $%
C^{\ast}(\Gamma;\Lambda _{T})$. Thus, the assertions in the following
theorem make sense.

\begin{theorem}
\label{CPident}The pair $(V,\pi )$ defined by equations (\ref{v1}) and (\ref
{v2}) is a faithful isometric covariant representation mapping $\mathcal{H}=%
\mathcal{H}_{T}$ into $C^{\ast }(\Gamma ;\Lambda _{T})$ in such a way that
its integrated form, $V\times \pi $, is a $C^{\ast }$-isomorphism mapping $%
\mathcal{O}_{\mathcal{H}_{T}}$ onto $C^{\ast }(\Gamma ;\Lambda _{T})$.
\end{theorem}

\begin{proof}
As we just mentioned, it is routine to check that $(V,\pi )$ is an isometric
covariant representation of $\mathcal{H}$ in $C^{\ast }(\Gamma ;\Lambda
_{T}) $. Also, it is evident that the image generates $C^{\ast }(\Gamma
;\Lambda _{T})$. Of course, $\pi $ is faithful and if $\beta :\mathbf{T}%
\rightarrow \mbox{\rm Aut}(C^{\ast }(\Gamma ;\Lambda _{T}))$ 
is the automorphism group induced by the position cocycle $c:\Gamma
\rightarrow \mathbf{Z}$, $c(x,m,y)=m$, then $\beta _{z}(f)(\lambda )=z^{c(%
\dot{\lambda})}f(\lambda )$, $f\in C_{c}(\Gamma ;\Lambda _{T})$, and it is
clear that $\beta _{z}\circ \pi =\pi $, while $\beta _{z}\circ V=zV$. Thus,
the only thing that needs to be checked is that $(V,\pi )$ satisfies the
Cuntz condition.

This, however, is quite easy. It is a matter of a couple of identifications
coupled with the appropriate references. First note that $K(\mathcal{H}_{T})=%
\mathcal{H}_{T}\otimes _{C_{0}(X)}\mathcal{H}_{T}^{\ast }$ which, in turn,
may be identified with the twisted groupoid $C^{\ast }$-algebra $C^{\ast
}(R(\sigma ),T\ast \overline{T}|_{R(\sigma )})$, where $R(\sigma
)=\{(x,y)\in X\times X\mid \sigma (x)=\sigma (y)\}$. This is a
straightforward computation, given the representation of $\mathcal{H}_{T}$
in terms of functions on $T$ that transform according to the formula: $\xi
(zt)=z\xi (t)$, $t\in T$, $z\in \mathbf{T}$. (See \cite{k3} and in
particular Section 3, therein, where a relation between the sheaf cohomology
of $X$ and the groupoid cohomology of $R(\sigma )$ is established.) Observe,
too, that $\pi ^{(1)}$ identifies $K(\mathcal{H}_{T})=C^{\ast }(R(\sigma
),T\ast \overline{T}|_{R(\sigma )})$ with the subalgebra of $C^{\ast
}(\Gamma ;\Lambda _{T})$ consisting of all elements that are supported on $%
\Lambda _{T}|_{X_{1,1}}$, in the notation of the proof of Theorem \ref
{twistiso}. Indeed, $\iota _{1,1}(X_{1,1})=R(\sigma )$.

Next, we claim that $\varphi (C_{0}(X))\subseteq K(\mathcal{H}_{T})$. This,
however, is obvious from the following facts: ($i$) $X$ may be identified as
the clopen subset $\Delta =\{(x,x)\in R(\sigma )\mid x\in X\}$; ($ii$) $%
T\ast \overline{T}|_{\Delta }$ is trivial (see the proof of Theorem \ref
{twistiso}); and ($iii$) in the identification of $K(\mathcal{H}_{T})$ with $%
C^{\ast }(R(\sigma ),T\ast \overline{T}|_{R(\sigma )})$, 
\begin{equation*}
\varphi (f)([t_{1},t_{2}])=\left\{ 
\begin{array}{cc}
f(p(t_{1}))\text{,} & [t_{1},t_{2}]\in T\ast \overline{T}|_{\Delta } \\ 
0 & \text{otherwise}
\end{array}
\text{.}\right.
\end{equation*}

Finally, we see that the Cuntz condition is satisfied by $(V,\pi )$ simply
by noting that the calculation of the previous paragraph and the
identification $\iota _{1,1}$ of $X_{1,1}$ with $R(\sigma )$ allows us to
identify $\pi ^{(1)}\circ \varphi $ with $\pi $.
\end{proof}

\bigskip

\section{Examples}

In this section we gather together several examples that illustrate some of
the theory we have developed.\bigskip

\noindent \textbf{Example 1} Let $X=\mathbf{T}^{k}$ be the $k$-dimensional
torus and $\sigma :X\rightarrow X$ be a covering map given by a $k$ by $k$
integer matrix $R$ with $|\det R\, |\geq 2$. Then $H^{0}(X,\mathbf{Z})=%
\mathbf{Z}$, and for $n\geq 1$, $H^{n}(X,\mathbf{Z})$ may be identified with 
$\mathbf{Z}^{k}\wedge \cdots\wedge \mathbf{Z}^{k}$, $n$ times. The map
induced on cohomology, $\sigma ^{n}$, can be identified with $R\wedge \cdots
\wedge R $, $n$ times. For $n=0$, the wedge product is taken to be the
identity map. From the exact sequence (\ref{2nd les}) we get 
\begin{equation*}
H^{0}(\Gamma ,\mathbf{Z})\simeq \mathbf{Z}
\end{equation*}
\begin{equation*}
H^{n}(\Gamma ,\mathbf{Z})\simeq \mbox{ker}\;(I-\underbrace{R\wedge \cdots
\wedge R}_{n})\oplus \mbox{coker}\;(I-\underbrace{R\wedge \cdots \wedge R}%
_{n-1}),\;n\geq 1.
\end{equation*}
Hence, by Corollary \ref{les bis}, 
\begin{equation*}
Br(\Gamma ) \simeq H^{3}(\Gamma ,\mathbf{Z})\simeq \mbox{ker}\;(I-R\wedge
R\wedge R)\oplus \mbox{coker}\;(I-R\wedge R).
\end{equation*}
In particular, the Brauer group can be infinite for $k\geq 3$.

\bigskip

\noindent \textbf{Example 2} Let $X$ be the infinite path space of the
topological graph 
\begin{equation*}
\mathbf{T}\overset{s}{\longleftarrow }\mathbf{T}\overset{r}{\longrightarrow }%
\mathbf{T},
\end{equation*}
where $s$ and $r$ are the covering maps given by $x\mapsto x^{p}$ and $%
x\mapsto x^{q}$, respectively. Then 
\begin{equation*}
X=\{(x_{1},x_{2},...)\in \mathbf{T}^{\mathbf{N}}\mid
(x_{n})^{q}=(x_{n+1})^{p},\ n\geq 1\}.
\end{equation*}
We assume that $|p|,|q|\geq 2$ and $(p,q)=1$. Then $X$ is a solenoid, 
\begin{equation*}
X=\lim_{\longleftarrow }(X_{m},\pi _{m})\text{,}
\end{equation*}
where $X_{m}$ is the space of paths of length $m$, and where the maps $\pi
_{m}:X_{m+1}\rightarrow X_{m}$ are the projections 
\begin{equation*}
\pi _{m}(x_{1},x_{2},...,x_{m+1})=(x_{1},x_{2},...,x_{m}).
\end{equation*}
The fact that $p$ and $q$ are relatively prime implies that $X$ is connected
(see [Br]). It is easy to see that each $X_{m}$ is homeomorphic to $\mathbf{T%
}$, and that the projections $\pi _{m}$ are given by the map $s$. Indeed,
the maps $f_{m}:\mathbf{T}\rightarrow X_{m}$, 
\begin{equation*}
f_{m}(x)=(s^{m}x,s^{m-1}rx,...,sr^{m-1}x,r^{m}x)
\end{equation*}
realize the homeomorphisms, and the diagram 
\begin{equation*}
\begin{array}{ccc}
\mathbf{T} & \overset{s}{\longrightarrow } & \mathbf{T} \\ 
f_{m+1}\downarrow \hspace{9mm} &  & \hspace{5.5mm}\downarrow f_{m} \\ 
X_{m+1} & \overset{\pi _{m}}{\longrightarrow } & X_{m}
\end{array}
\end{equation*}
is commutative.

Moreover, in this identification, the unilateral shift $\sigma:X\rightarrow
X $ is given by

\begin{equation*}
\begin{array}{ccccccccc}
\mathbf{T} & \overset{s}{\longleftarrow} & \mathbf{T} & \overset{s}{%
\longleftarrow} & \mathbf{T} & \overset{s}{\longleftarrow} & ... & {%
\longleftarrow} & X \\ 
{} & {} & {} & {} & {} & {} & {} & {} & {} \\ 
{} & r\swarrow & {} & r\swarrow & {} & {} & {} & {} & \sigma\downarrow \;\;\;
\\ 
{} & {} & {} & {} & {} & {} & {} & {} & {} \\ 
\mathbf{T} & \overset{s}{\longleftarrow} & \mathbf{T} & \overset{s}{%
\longleftarrow} & \mathbf{T} & \overset{s}{\longleftarrow} & ... & {%
\longleftarrow} & X.
\end{array}
\end{equation*}

This allows us to calculate the integer cohomology of $X$ and to identify
the maps induced by $\sigma$: 
\begin{equation*}
H^{0}(X,\mathbf{Z})=\lim_{\longrightarrow}(\mathbf{Z},id)=\mathbf{Z}%
;\;\;\sigma^{0}=id,
\end{equation*}
\begin{equation*}
H^{1}(X,\mathbf{Z})=\lim_{\longrightarrow}(\mathbf{Z},p)=\mathbf{Z}[
1/p];\;\; \sigma^{1}=q/p,
\end{equation*}
\begin{equation*}
H^{k}(X,\mathbf{Z})=0, \ k\geq 2.
\end{equation*}

From the long exact sequence (\ref{2nd les}) we determine the integer
cohomology of the groupoid $\Gamma =\Gamma (X,\sigma )$: 
\begin{equation*}
H^{0}(\Gamma ,\mathbf{Z})=\mathbf{Z},
\end{equation*}
\begin{equation*}
H^{1}(\Gamma ,\mathbf{Z})=\mathbf{Z}\oplus \mbox{ker}\;(1-q/p)=\mathbf{Z},
\end{equation*}
\begin{equation*}
H^{2}(\Gamma ,\mathbf{Z})=\mbox{coker}\;(1-q/p)=\mathbf{Z}/(p-q)\mathbf{Z},
\end{equation*}
\begin{equation*}
H^{k}(\Gamma ,\mathbf{Z})=0, \ k\geq 3.
\end{equation*}

In particular, by Corollary \ref{les bis} again, 
\begin{equation*}
Br(\Gamma ) \simeq H^{3}(\Gamma ,\mathbf{Z})=0.
\end{equation*}
\bigskip

\noindent\textbf{Example 3} Given a sequence of local homeomorphisms as in 
\cite[Addendum 3]{k2} 
\begin{equation*}
X_0 \overset{\sigma_0}{\longrightarrow} X_1 \overset{\sigma_1}{%
\longrightarrow} X_2 \overset{\sigma_2}{\longrightarrow} \cdots,
\end{equation*}
take X to be the disjoint union of the spaces, $X = \coprod_k X_k$ and
define $\sigma : X \rightarrow X$ in the natural way: if $x \in X_k \subset
X $ set $\sigma (x) = \sigma_k(x)$. Let $\Gamma = \Gamma(X, \sigma)$; if $%
\sigma_n$ is surjective for all $n$, then $X_0$ meets every orbit. It
follows that the reduction $\Gamma|_{X_0}$ is equivalent to $\Gamma$ and
therefore has the same cohomology. Further, $\Gamma|_{X_0}$ is precisely the
ultraliminary groupoid considered in \cite[Addendum 3]{k2} (the equivalence
relation on $X_0$ induced by the maps $\sigma_n\cdots\sigma_0$). We show how
Corollary \ref{les} allows one to recover the short exact sequence for the
cohomology given in \cite[Addendum 3]{k2}.

A $\Gamma$-sheaf $A$ is given by a sequence of sheaves $A_k$ over $X_k$
together with identifications $A_k = \sigma_k^*(A_{k+1})$ (see the proof of
Corollary \ref{les}). Given such a $\Gamma$-sheaf $A$, we have $H^n(X, A) =
\prod_k H^n(X_k, A_k)$ and $\sigma^* = \prod_k \sigma_k^*$ (the map induced
on cohomology). By Corollary \ref{les}, $H^{0}(\Gamma, A)$ is isomorphic to
the kernel of the map 
\begin{equation*}
1 - \prod_k \sigma_k^* : \prod_k H^0(X_k, A_k) \rightarrow \prod_k H^0(X_k,
A_k);
\end{equation*}
that is, $H^{0}(\Gamma, A)$ is isomorphic to the subgroup consisting of all $%
(g_k) \in \prod_k H^0(X_k, A_k)$ for which $g_k = \sigma_k^*(g_{k+1})$.
Hence, 
\begin{equation*}
H^{0}(\Gamma, A) = \lim_{\longleftarrow} H^0(X_k, A_k).
\end{equation*}
Similarly, for $n > 0$, $H^{n}(\Gamma, A)$ is an extension of the cokernel
of the map 
\begin{equation*}
1 - \prod_k \sigma_k^* : \prod_k H^{n-1}(X_k, A_k) \rightarrow \prod_k
H^{n-1}(X_k, A_k)
\end{equation*}
by the kernel of the map 
\begin{equation*}
1 - \prod_k \sigma_k^* : \prod_k H^n(X_k, A_k) \rightarrow \prod_k H^n(X_k,
A_k).
\end{equation*}
Hence, we obtain the short exact sequence 
\begin{equation*}
0 \rightarrow {\lim_{\longleftarrow}}^1 H^{n-1}(X_k, A_k) \rightarrow
H^{n}(\Gamma, A) \rightarrow \lim_{\longleftarrow} H^n(X_k, A_k) \rightarrow
0.
\end{equation*}

\bigskip

\noindent \textbf{Example 4} (Skew product construction) The following
construction is adapted from \cite{kp}. Given a local homeomorphism $\sigma
:X\rightarrow X$ and a continuous map $c:X\rightarrow G$, where $G$ is a
locally compact group, one constructs a new local homeomorphism $\tau
:X\times G\rightarrow X\times G$ by the formula 
\begin{equation*}
\tau (x,g)=(\sigma (x),gc(x)).
\end{equation*}
We define a continuous one-cocycle $\tilde{c}:\Gamma (X,\sigma )\rightarrow
G $ in the following way (see lemma \ref{naturaliso}). For $\gamma
=(x,k-l,y)\in \Gamma $ with $\sigma ^{k}(x)=\sigma ^{l}(y)$, set 
\begin{equation*}
\tilde{c}(\gamma )=c(x)c(\sigma (x))\cdots c(\sigma ^{k-1}(x))c(\sigma
^{l-1}(y))^{-1}\cdots c(\sigma (y))^{-1}c(y)^{-1}
\end{equation*}
($\tilde{c}$ is clearly well-defined and satisfies the cocycle property).
One may construct the skew-product as defined by Renault, $\Gamma (X,\sigma
)(\tilde{c})=\Gamma (X,\sigma )\times _{\tilde{c}}G$ (see \cite[Def I.1.6]
{re}); it is straightforward to verify that 
\begin{equation*}
\Gamma (X\times G,\tau )\simeq \Gamma (X,\sigma )\times _{\tilde{c}}G,
\end{equation*}
see \cite[2.4]{kp}. It follows by \cite[II.5.7]{re} that when $G$ is abelian 
\begin{equation*}
C^{\ast }(\Gamma (X\times G,\tau ))\simeq C^{\ast }(\Gamma (X,\sigma
))\times _{\alpha _{\tilde{c}}}\hat{G},
\end{equation*}
where $\alpha =\alpha _{\tilde{c}}$ is the action of $\hat{G}$ on $C^{\ast
}(\Gamma (X,\sigma ))$ induced by the cocycle. In particular, taking $G=%
\mathbf{R}$, $X=\{1,2,...,n\}^{\mathbf{N}}$, $\sigma $ the Bernoulli shift,
and $\lambda _{1},\dots ,\lambda _{n}\in \mathbf{R}$, we may define the
continuous function $c:X\rightarrow \mathbf{R}$ by $c(x)=\lambda _{k}$ if $%
x_{1}=k$. Then $C^{\ast }(\Gamma (X,\sigma ))\simeq \mathcal{O}_{n}$ and the
induced action of $\mathbf{R}$ on $\mathcal{O}_{n}$ is given by $\alpha
_{t}(S_{k})=e^{it\lambda _{k}}S_{k}$; the associated crossed product $%
\mathcal{O}_{n}\times \mathbf{R}$ is a special case of those studied by
Kishimoto in \cite[\S 4]{ki} (see also \cite{kk}). 

\end{document}